\documentclass[11pt]{article}
\usepackage[left=1in,right=1in,top=1in,bottom=1in]{geometry}
\usepackage{times}
\usepackage{expl3}
\usepackage{cite}
\usepackage[table]{xcolor}
\usepackage{multirow}
\usepackage{stackengine} 
\usepackage{hhline}
\usepackage{lipsum}
\usepackage{titlesec}
\usepackage{wrapfig}
\usepackage{epsfig}
\usepackage{graphicx}
\usepackage{amsmath}
\usepackage[title]{appendix}
\usepackage{amssymb}
\usepackage{epstopdf}
\usepackage{boldline}
\usepackage{calligra}
\usepackage{url}
\usepackage{blindtext}

\newcommand{\define}{\stackrel{\mbox{\tiny def}}{=}}

\newtheorem{definition}{Definition}
\newtheorem{theorem}{Theorem}

\newtheorem{lemma}{Lemma}

\usepackage{mathtools}
\usepackage{epstopdf}
\usepackage{balance}
\usepackage{thmtools}
\usepackage{thm-restate}
\usepackage{hyperref}
\usepackage{cleveref}

\usepackage[ruled,vlined]{algorithm2e}
\include{pythonlisting}
\newcommand{\ostar}{\mathbin{\mathpalette\make@circled\star}}

\makeatletter
\newcommand{\removelatexerror}{\let\@latex@error\@gobble}
\makeatother
\makeatletter
\newcommand*{\rom}[1]{\expandafter\@slowromancap\romannumeral #1@}
\makeatother

\ExplSyntaxOn
\newcommand\latinabbrev[1]{
  \peek_meaning:NTF . {% Same as \@ifnextchar
    #1\@}%
  { \peek_catcode:NTF a {% Check whether next char has same catcode as \'a, i.e., is a letter
      #1.\@ }%
    {#1.\@}}}
\ExplSyntaxOff

\titleclass{\subsubsubsection}{straight}[\subsubsection]
\date{}
\begin{document}
\vspace{1cm}
\title{Geometric Approach For Majorizing Measures on Hadamard Manifolds}\vspace{1.8cm}
\author{Shih~Yu~Chang 
% <-this % stops a space
\thanks{Shih Yu Chang is with the Department of Applied Data Science,
San Jose State University, San Jose, CA, U. S. A. (e-mail: {\tt
shihyu.chang@sjsu.edu}).
           }}

\maketitle

\begin{abstract}
Gaussian processes can be treated as subsets of a standard Hilbert space, however, the volume size relation between the underlying index space of random processes and its convex hull is not clear. The understanding of such volume size relations can help us to establish a majorizing measure theorem geometrically. In this paper, we assume that the underlying index space of random processes is a simply connected manifold with sectional curvature less than negative one (Hadamard manifold). We derive the upper bound for the ratio between the volume of the underlying index space and the volume of its convex hull. We then apply this volume ratio to prove the majorizing measure theorem geometrically. 
\end{abstract}


\section{Introduction}\label{sec:Introduction} 

Majorizing measures provide estimates for the supremum of stochastic processes. The chaining argument used in majorizing measures can be traced back to Kolmogorov. Although majorizing measures stem from the theory of Gaussian processes at the beginning times, nowaday, these measures have been extended beyond this Gaussian process assumption. Let us review the following definition about majorizing measures.
\begin{definition}\label{def:gamma 2}
Given a metric space $(T, d(\cdot, \cdot))$, we define $\gamma_{\alpha}(T,d(\cdot, \cdot))$ as  
\begin{eqnarray}
\gamma_\alpha(T, d(\cdot, \cdot)) =\inf \sup\limits_{t \in T} \sum\limits_{m \geq 0}2^{m/\alpha} \Delta(\mathrm{A}_n(t)),  
\end{eqnarray}
where the infimum is taken over all admissible sequences and $\Delta(\mathrm{A}_n(t))$ denotes the diameter of $\mathrm{A}_n(t)$ with respect to the metric $d(\cdot, \cdot)$. 
\end{definition}
In the Gaussian case, that is when $X_t = \sum\limits_{t_k \in T} t_k g_k$, where $g_k$ are i.i.d. standard Gaussians, we have the celebrated Fernique–Talagrand majorizing measure theorem as:
\begin{eqnarray}
\frac{1}{L} \gamma_2(T, d(\cdot, \cdot)) \leq \mathbb{E} \sup\limits_{t \in T} X_t   \leq L \gamma_2(T, d(\cdot, \cdot)).
\end{eqnarray}
The majorizing measure theorem is a main method used to prove Theorem 2.11.1 in~\cite{talagrand2021upper}. 
Gaussian processes can be seen as subsets of a standard Hilbert space, but the
geometric understanding that would relate the size of a set with the size of its
convex hull is still insufficient~\cite{talagrand2021upper}. The purpose of this work is 
to establish Theorem 2.11.1 in~\cite{talagrand2021upper} based on a geometrical approach, i.e., we wish to show
\begin{eqnarray}\label{eq:2.130}
\gamma_2 (T_h) \leq L \gamma_2(T),
\end{eqnarray}
where $T_h$ is the convex hull for the original space $T$. 

In this work, we assume that the space $T$ is a simply connected manifold wih sectional curvature $k < -1$, and we use $T_h$ to represent the convex hull of the space $T$. We first use the method adopted by~\cite{borbely1998some} to provide a more general upper bound estimation for the volume of $T_h$. The next step is to give a lower bound for the volume estimation of the space $T$ according to Theorem 3 in~\cite{borisenko2002convex}. Then, the covering number ratio between the space $T$ and the space $T_h$ can be obtained through the volume ratio between the space $T$ and the space $T_h$. This ratio will help us to establish Theorem 2.11.1 in~\cite{talagrand2021upper}. 

The rest of this paper is organized as follows. The upper bound estimation for the volume $T_h$ and the lower bound estimation for the volume $T$ are given in Section~\ref{sec:Volume Estimation for T and Th}. We bound the covering number of $T_h$ in terms of the covering number of $T$ in Section~\ref{sec:Covering Number Bounds by Volume}. Finally, the majorizing measuring theorem is proved from a geometric perspective in Section~\ref{sec:Geometric Proof}.

\section{Volume Estimation for Spaces $T$ and $T_h$}\label{sec:Volume Estimation for T and Th}

\subsection{Volume Estimation of $T_h$}\label{sec:Volume of Th}

The purpose of this section is to derive the volume estimate for $T_h$. Our approach is based on the method provided by~\cite{borbely1998some}, but we fix some issues and relax some assumptions from there. 

We begin to restate Lemma 2 and Proposition 1 from~\cite{borbely1998some} as the following two Lemmas for later proof presentation convenience. 

\begin{lemma}\label{lma:Lemma2}
Let $M$ be a Hadamard manifold with dimension $n$ and sectional curvatures $K$ within in the range
$-k^2 \leq K \leq -1$. We also have $V_1,\cdots,V_m$ are convex 
sets in $M$ with $\bigcap_{i=1}^{m} V_i \neq \empty$ and $T = \bigcup_{i=1}^{m} V_i$. Then, there is a 
constant $C(k)$ depending only on the pinching, such that the convex hull of $T$ lies in the C-neighborhood of $T$.
\end{lemma}

\begin{lemma}\label{lma:Prop1}
Let $M$ be a Hadamard manifold with dimension $n$ and sectional curvatures $K$ within in the range
$-k^2 \leq K \leq -1$. Also let $\mathbf{B}_r$ be a closed geodesic ball with radius $r >0$ and $G \in \mathbf{B}_r$ be a closed convex set. Then, there is a contant $C(n,k,r)$, depending on the dimension of $M$, the sectional curvature, and the radius $r$, such that 
\begin{eqnarray}
\mbox{Vol}(N(G, \delta) - G) < \delta C(n,k,r),
\end{eqnarray}
where $N(G, \delta)$ is the $\delta$-neighborhood of $G$ in $\mathbf{B}_r$
\end{lemma}

\begin{theorem}\label{thm:upper bound of vol Th}
Let $M$ be a Hadamard manifold with dimension $n$ and sectional curvatures $K$ within in the range
$-k_2^2 \leq K \leq -k_1^2$, where $k_2 >  k_1 > 1$. We also have $V_1,\cdots,V_m$ are $m$ $\lambda$-convex 
sets in $M$ with $\lambda \leq k_2$ and $\bigcap_{i=1}^m V_i \neq \emptyset$. Suppose we select $m$ points, $P_1, P_2, \cdots, P_m$, such that all these $P_i$ for $1 \leq i \leq m$ are sampled from the space $T \define \bigcup_{i=1}^{m^\rho} V_i$, where $0 \leq \rho \leq 1$. Let us define $T_h$ as the convex hull of space $T$, then we have  
\begin{eqnarray}\label{eq1:thm:upper bound of vol Th}
\mbox{Vol}(T_h) \leq C_{ub} m^{1+\rho- \varpi(\rho, k_1,k_2,n)}.
\end{eqnarray}
where $\varpi(\rho, k_1,k_2, n)$ is a constant depending on the sectional curvature bounds and the manifold dimension. 
\end{theorem}

\textbf{Proof:}

The proof of this Theorem is composed of two portions. The first portion is to construct the convex set to cover $T$, and the second portion is to provide a volume upper bound estimate of the constructed convex set. 

The proof is based on induction and the crucial step is to understand how the volume is changed from a convex hull by adding one more point. Without loss of generality, we assume that $P_1, \cdots, P_m \in S_{\infty}(M)$,    where $S_{\infty}(M)$ denotes the ideal boundary. The convex hull made by points $P_1, \cdots, P_{m-1}$ is denoted by $\tilde{T}_h$, and the convex hull made by the space $\tilde{T}_h$ with one more point $P_m$ is denoted by $T_h$. 

Let $Q \in \partial \tilde{T}_h$ (boundary of the space $\tilde{T}_h$) be the closest point to $P_m$, and use $\zeta(t)$ for $t \in [0, \infty)$to represent the unit speed geodesic ray connecting $Q$ and $P_m$. Note that the line $\zeta(t)$ is perpendicular to the boundary $\partial \tilde{T}_h$. We say that $\partial \tilde{T}_h$ is not smooth if the angle $\angle R Q \Upsilon$ satisfies $\angle R Q \Upsilon \geq \frac{\pi}{2}$ for every $R \in \zeta(t)$ and $\Upsilon \in \tilde{T}_h$.

We will construct a convex set $G$ such that $T_h \subset G$. We will define the following two functions $g_1(P)$ and $g_2(P)$ for $P \in M$ first and apply them to construct the space $G$. The function $g_1(P)$ is defined as
\begin{eqnarray}\label{eq:g1 def}
g_1(P) = 1 - \exp( - a \times \mbox{dist}(P, \tilde{T}_h)),
\end{eqnarray}
where $a$ is a positive constant. Moreover, we define the function $g_2(P)$ as
\begin{eqnarray}\label{eq:h2 def}
g_2(P) = 1 - \exp( - a \times \mbox{dist}(P, \zeta(t))).
\end{eqnarray}
Then, we can define the space $G$ as
\begin{eqnarray}\label{eq:H construction}
G = \{P \in M: g_1(P) + g_2(P) \leq 1 \}. 
\end{eqnarray}

Since $\partial G$ is not necessarily smooth, we cannot show the convexity of $G$ by the second fundamental form  positive definite. Therefore, we can construct a smooth enveloping surface with positive definite second fundamental form that includes $G$. For every point $P \in \partial G$, we will construct a space $\acute{G} \supset G$ such that $P \in \partial \acute{G}$, where $\partial \acute{G}$ is smooth near $P$ and $\partial \acute{G}$ has a positive second fundamental form at $P$. This will show the convexity of $G$. 

\begin{figure}[htbp]
	\centerline{\includegraphics[width=0.9\columnwidth,draft=false]
		{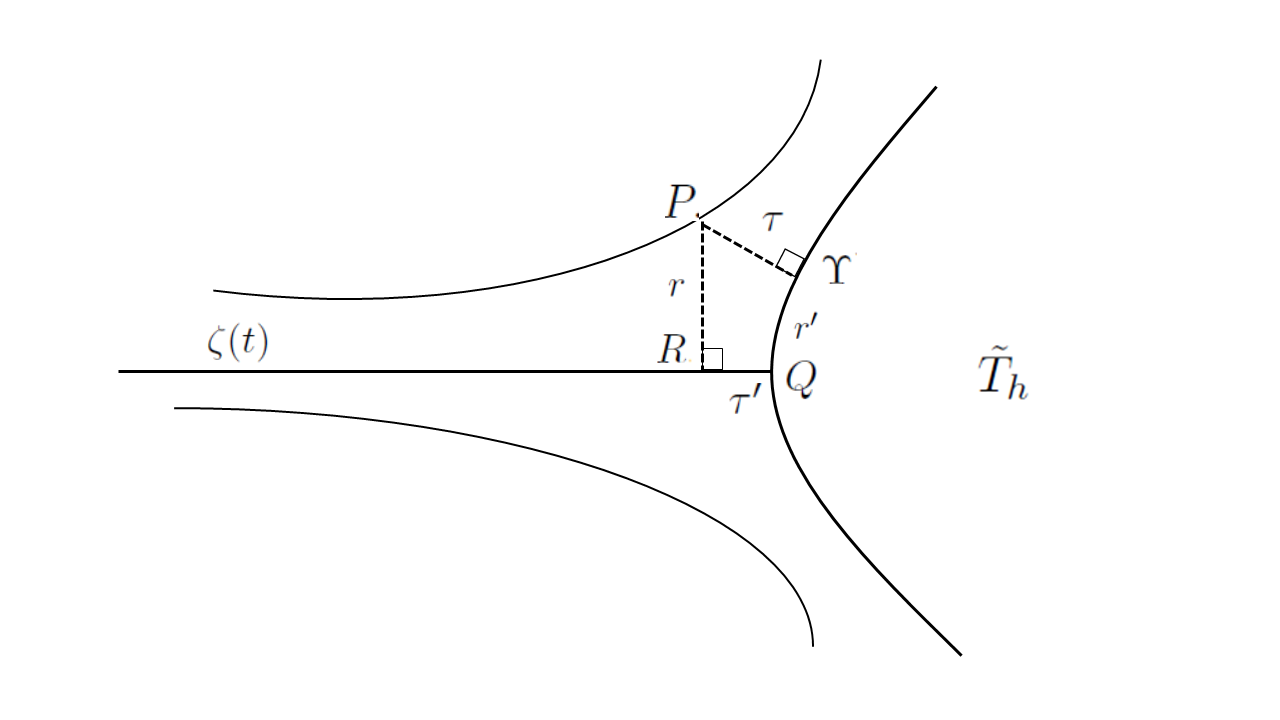}}
	\caption{Proof plot.}\label{fig:ProofPlot}
\end{figure}

Let $P \in \partial G$, and $R$ and $\Upsilon$ be two closest points to $P$ on $\zeta(t)$ and $\tilde{T}_h$, respectively. We use $r$ to represent $\mbox{dist}(P, R)$, and $\tau$ to represent $\mbox{dist}(P, \Upsilon)$, see Fig.~\ref{fig:ProofPlot}. Let $V_1$ be the closed half-space at $\Upsilon$ containing $\tilde{T}_h$ obtained by collecting geodesic rays starting from $\Upsilon$ and having angle greater or equal to $\frac{\pi}{2}$ with the geodesic segment $\Upsilon P$. If we define the function $\acute{g}_1$ as
\begin{eqnarray}\label{eq:acute g1 def}
\acute{g}_1 \define 1 - \exp( - a \times \mbox{dist}(P, V_1)).
\end{eqnarray}
Then, we have a smooth surface $\partial V_1$ such that $V_1 \supset \tilde{T}_h$ and $\acute{g}_1 \leq g_1$. Similarly, let $V_2$ be the closed half-space at $Q$ perpendicular to $\zeta(t)$ containing $\tilde{T}_h$. If $P \notin \partial V_2$, then $h_2$ is smooth near $P$ and we set $\acute{g}_2 = g_2$.  If $P \in \partial V_2$, then let $\acute{\zeta}(t)$ be the extension of geodesic ray $\zeta(t)$ beyond $Q$; that is, for some small $\epsilon > 0$, $\acute{\zeta}(t) = \zeta[-\epsilon, \infty]$ and set the function $\acute{g}_2$ as
\begin{eqnarray}\label{eq:acute g2 def}
\acute{g}_2 \define 1 - \exp( - a \times \mbox{dist}(P, \acute{\zeta}(t))).
\end{eqnarray}
Then, we have $\acute{\zeta}(t) \supset \zeta(t)$ and $\acute{g}_2 \leq g_2$.

We also define $\acute{G} = \{R \in M: \acute{g}_1(R) + \acute{g}_2(R) \leq 1 \}$. Since $\acute{g}_1(R) + \acute{g}_2(R) \leq g_1(R) + g_2(R) $ for any $ R (\neq P) \in M$ and $\acute{g}_1(P) + \acute{g}_2(P) =  g_1(P) + g_2(P)$, we have $\acute{G} \supset G$ and $P \in \partial \acute{G} \cap \partial G$. Our next goal is to show that $\partial \acute{G}$ has positive definite second fundamental form at $P$.

From the above construction, we have $\angle P \Upsilon Q \geq \frac{\pi}{2}$ and $\angle PRQ = \frac{\pi}{2}$, therefore, we have $\angle RP \Upsilon \leq \frac{\pi}{2}$ from the negative curvature assumption of $M$. If we define the following two normal vectors $N_1$ and $N_2$ as
\begin{eqnarray}
N_1 \define \frac{\nabla \mbox{dist}(P, V_1)          }{\left\Vert  \nabla \mbox{dist}(P, V_1)     \right\Vert},
\end{eqnarray}
and 
\begin{eqnarray}
N_2 \define \frac{\nabla \mbox{dist}(P, \acute{\zeta}(t))          }{\left\Vert  \nabla \mbox{dist}(P, \acute{\zeta}(t) )     \right\Vert};
\end{eqnarray}
then, we have
\begin{eqnarray}\label{eq:21}
\langle N_1, N_2 \rangle > 0,
\end{eqnarray}
due to the angle $\angle RP \Upsilon \leq \frac{\pi}{2}$. 

Let $X$ be a unit tangent vector to $\partial \acute{G}$ at $P$, then we have 
\begin{eqnarray}\label{eq:22}
\langle e^{-a \tau} N_1 + e^{-a r} N_2, X \rangle = 0,
\end{eqnarray}
where $e^{-a \tau} + e^{-a r} = 1$. From Eq.~\eqref{eq:21}, we have
\begin{eqnarray}
\left\Vert e^{-a r} N_1 - e^{-a \tau} N_2 \right\Vert \leq \sqrt{e^{-2 a r} + e^{-2 a \tau} }.
\end{eqnarray}
Then, we have the following two inequalities:
\begin{eqnarray}\label{eq:23-a1}
\left\vert \langle N_1, X \rangle   \right\vert \leq \frac{e^{- a r }}{ \sqrt{e^{-2 a r} + e^{-2 a \tau} } },
\end{eqnarray}
and
\begin{eqnarray}\label{eq:23-a2}
\left\vert \langle N_2, X \rangle   \right\vert \leq \frac{e^{- a \tau}}{ \sqrt{e^{-2 a r} + e^{-2 a \tau} } }.
\end{eqnarray}
Also, since $g_1(P) + g_2(P) = 1$, we also have $e^{- a r} + e^{- a \tau} = 1$. Then, from Eqs.~\eqref{eq:23-a1} and~\eqref{eq:23-a2}, we have 
\begin{eqnarray}\label{eq:23-b1}
\left\vert \langle N_1, X \rangle   \right\vert \leq \sqrt{2} e^{-a r},
\end{eqnarray}
and
\begin{eqnarray}\label{eq:23-b2}
\left\vert \langle N_2, X \rangle   \right\vert \leq \sqrt{2} e^{-a \tau}.
\end{eqnarray}

The next thing required for the convexity proof is to estimate the second derivative of $\acute{g} \define = \acute{g}_1 + \acute{g}_2$ at $P$ in the direction of $X$. Since we have the following relation for the second 
derivative of the function $\acute{g}(P)$ at the direction $X$, denoted as $D_{X}^2 \acute{g}(P)$,  
\begin{eqnarray}\label{eq:24-2}
D_{X}^2 \acute{g}(P) &=& -a^2 e^{-a \tau}\langle N_1, X \rangle^2 + 
ae^{-a \tau}D_{X}^2 \mbox{dist}(P, V_1)\nonumber \\
&  &  - a^2 e^{-a r}\langle N_2, X \rangle^2 + 
ae^{-a r}D_{X}^2 \mbox{dist}(P, \acute{\zeta}(t) ).
\end{eqnarray}

Since the sectional curvature is upper bounded by $k_1^2$, by applying comparison technique we have
\begin{eqnarray}\label{eq:24-1a}
D_{X}^2 \mbox{dist}(P, V_1) \geq (k_1^2 - \langle N_1, X \rangle^2) \tanh(\mbox{dist}(P, V_1)),
\end{eqnarray}
and
\begin{eqnarray}\label{eq:24-1b}
D_{X}^2 \mbox{dist}(P, \acute{\zeta}(t)) \geq (k_1^2 - \langle N_2, X \rangle^2) \tanh( \mbox{dist}(P, \acute{\zeta}(t)) ).
\end{eqnarray}
From Eqs.~\eqref{eq:24-1a} and~\eqref{eq:24-1b}, Eq.~\eqref{eq:24-2} can be further expressed as
\begin{eqnarray}
D_{X}^2 \acute{g}(P) &\geq& -a^2 ( e^{-a \tau}\langle N_1, X \rangle^2 + 
 e^{-a r }\langle N_2, X \rangle^2) + a e^{-a \tau}( k_1^2 - \langle N_1, X \rangle^2 ) \tanh(\tau)
\nonumber \\
&  &+ a e^{- a r } ( k_1^2 - \langle N_2, X \rangle^2 ) \tanh(r).
\end{eqnarray}

In order to show the positive definite for the Hessian of $\acute{G}$, we have to show 
\begin{eqnarray}\label{eq:24+1}
a e^{-a \tau}( k_1^2 - \langle N_1, X \rangle^2 ) \tanh(\tau)
+ a e^{- a r } ( k_1^2 - \langle N_2, X \rangle^2 ) \tanh(r) >  \nonumber \\
a^2 ( e^{-a \tau}\langle N_1, X \rangle^2 + 
 e^{-a r }\langle N_2, X \rangle^2).
\end{eqnarray}
From Eqs.~\eqref{eq:23-b1},~\eqref{eq:23-b2} and~\eqref{eq:24+1}, the positive definite for the Hessian of $\acute{G}$ can be proved if we have
\begin{eqnarray}\label{eq:24+2}
2a^2 e^{-a (\tau + r)} < a e^{-a \tau}( k_1^2 - \langle N_1, X \rangle^2 ) \tanh(\tau)
+ a e^{- a r } ( k_1^2 - \langle N_2, X \rangle^2 ) \tanh(r).
\end{eqnarray}

Since Eq.~\eqref{eq:24+2} is symmetric with respect to $r$ and $\tau$, and we have $e^{-a r} + e^{-a \tau} = 1$, we can assume that $e^{- a \tau} < \frac{1}{2}$. Moreover, both terms in Eq.~\eqref{eq:24+2} are positive due to that $N_1, N_2$ and $X$ are unit vectors. Therefore, the convexity of the set $G$ is equivalent to show that we can find the proper positive constant $a$ to satisfy the following expression:
\begin{eqnarray}\label{eq:24+3}
2a^2 e^{-a (\tau + r)} < a e^{-a \tau}( k_1^2 - \langle N_1, X \rangle^2 ) \tanh(\tau).
\end{eqnarray}
From Eqs.~\eqref{eq:23-b1},~\eqref{eq:23-b2} and $e^{- a r} \leq \frac{1}{2}$, we can re-write Eq.~\eqref{eq:24+3} as
\begin{eqnarray}\label{eq:24+4}
2a^2 e^{-a (\tau + r)} < a e^{-a \tau}( k_1^2 - \frac{1}{2} ) \tanh(\tau).
\end{eqnarray}

By specifying the value range of $\tau$, we can have following two cases: (I) $\tanh(\tau) \geq \frac{\tau}{2}$, and  (II) $\tanh(\tau) \geq \frac{1}{2}$. 

For the case (I), we have $e^{- a r } = 1 - e^{-a \tau} < a \tau$, then Eq.~\eqref{eq:24+4} becomes  
\begin{eqnarray}\label{eq:25a}
2a^2 \tau < ( k_1^2 - \frac{1}{2}) \tanh(\tau).
\end{eqnarray}
For the case (II), Eq.~\eqref{eq:24+4} becomes  
\begin{eqnarray}\label{eq:25b}
\frac{2ae^{-a r}}{ k_1^2 - \frac{1}{2} } < \frac{1}{2}.
\end{eqnarray}
By setting the value of $k_1 \geq 1$, we can find the value $a$ to satisfy both Eqs.~\eqref{eq:25a} and~\eqref{eq:25b} . For example, if we set $k_1 = 1$, we can find $a = \frac{1}{4}$ to satisy both Eqs.~\eqref{eq:25a} and~\eqref{eq:25b}. Therefore, we finish the proof of the first portion about the construction of the convex set $G$ to cover $T_h$.

Now, we begin the second portion of the proof about the volume upper bound estimate of the constructed convex set $G$. 

We note that there is an constant $\eta_0 > 0$, such that, for every $\eta > \eta_0$, the ball with the radius $\eta$ centering at $Q$, denoted as $B(Q,\eta)$, disconnects the set $G$ such that no point of the ray $\zeta - B(Q, \eta)$ can be connected to any point of $\tilde{T}_h - B(Q, \eta)$ by a curve lying entirely in the set $G - B(Q,\eta)$. The estimation of $\mbox{Vol}(G)$ will be separated into three parts. The first part, denoted as $G_1$, is the connected component of $G - B(Q,\eta)$ containing $\zeta - B(Q,\eta)$ (the long thin tube around $\zeta(t)$, see Fig.1). The second part, denoted as $G_2$, is defined as $G_2 = G \bigcap B(Q,\eta)$. Finally, the third part, denoted as $G_3$, is the rest of $G - B(Q,\eta)$. 

We begin with the volume estimation of $G_1$. For $P \in \partial G_1 - B(Q,\eta)$, let $R$ and $\Upsilon$ as closest points to $P$ on the ray $\zeta$ and $\tilde{T}_h$, respectively. We also adopt the following distance notations: $r \define \mbox{dist}(P, R)$,  $r' \define \mbox{dist}(Q, \Upsilon)$,  $\tau \define \mbox{dist}(P, \Upsilon)$,  $\tau' \define \mbox{dist}(R, Q)$, see Fig. 1. Because the curvature is negative, we have $r' < r$ and $\tau' < \tau$. From the construction of $G_1$ and the condition of $g_1(P) + g_2(P) = 1$, we have 
\begin{eqnarray}
r &<& \tau, \nonumber \\
\eta &\leq & \mbox{dist}(R, Q) \leq r' + \tau < r + \tau, \nonumber \\
e^{-a r} + e^{-a \tau} &=& 1.
\end{eqnarray}
By selecting large enough $\eta_0$, we have 
\begin{eqnarray}\label{eq:26}
r < \frac{e^{-a(\tau - 1)}}{a} < \frac{e^{-a(\tau' - 1)}}{a}.
\end{eqnarray}

By volume comparison with the ray $\zeta$, we can estimate the volume of the tube $G_1$ with the radius bounded by Eq.~\eqref{eq:26} as
\begin{eqnarray}\label{eq:27}
\mbox{Vol}(G_1) < C_1(n, k_2) \int_{\eta-2}^{ \infty } \left( \frac{e^{-a(\tau' -1 )}}{a} \right)^{n-1} d \tau' < C_2(a, n, k_2) e^{-a(n-1)\eta},
\end{eqnarray}
where $C_1(n, k_2), C_2(a, n, k_2)$ are appropriate constants depending on the dimension
and the curvature of the underlying Hadamard manifold. Note that in Eq. (2.7) from~\cite{borbely1998some}, the Hadamard manifold dimension should be considered at the integrand.   
 
Let us estimate the volume of $G_2 = G \bigcap B(Q,\eta)$. It is easy to obtain since we have 
\begin{eqnarray}\label{eq:210}
\mbox{Vol}(G_2) < \mbox{Vol}(B(Q,\eta)) < C_3(n, k_2)e^{k_2 (m-1) \eta}.
\end{eqnarray}

We will beging to estimate the volume of $G_3$, which is the most complicated part in $G$ to estimate its volume. For induction argument, we will consider $\mbox{Vol}(G_3 - F)$. From the construction of $G_3$, we have
\begin{eqnarray}\label{eq:28}
\tau < \frac{e^{-a(r - 1)}}{a} < \frac{e^{-a(r' - 1)}}{a}.
\end{eqnarray}
Then, the thickness of $G_3 - \tilde{T}_h$ approaches to zero exponentially as $r'$ approaches to infinity.  To estimate the $\mbox{Vol}(G_3 - F)$, we cut $(G_3 - F)$ into bounded pieces. According to Lemma~\ref{lma:Lemma2}, we can find the constant $C(k)$ such that $ \tilde{T}_h$ lies in the tubular neighbor of the rays $\overline{Q P_i}$, where $i=1,2,\cdots, m-1$. Moreover, from the construction of $G$, we know that $G$ remains within a distance $4 \ln 2 $ of the set $\zeta \bigcup \tilde{T}_h$. Let us use $Q_i(x)$ to represent the point on the geodesic ray $\overline{Q P_i}$ for $i=1,2,\cdots, m-1$ such that $\mbox{dist}(Q, Q_i(x)) = x$. Then, if $\eta_0$ is larger than a sufficiently large absolute constant and pick a positive value $\theta > 4 \ln 2$, the set $\tilde{T}_h - B(Q, \eta)$ and $G_3$ will be covered by the collection of balls $B(Q_i(x_{\ell}), 2 \theta)$, where $i=1,2,\cdots, m$ and $x_{\ell} = \eta + \frac{\ell \theta}{2}$ for $\ell=1,2,\cdots$. We will estimate $\mbox{Vol}( (G_3 - \tilde{T}_h) \bigcap B(Q_i(x_{\ell}), 2 \theta) )$ with respect to each ball $B(Q_i(x_{\ell}), 2 \theta) )$, individually. 

From Lemma~\ref{lma:Prop1} and Eq.~\eqref{eq:28}, we have 
\begin{eqnarray}
\mbox{Vol}( (G_3 - \tilde{T}_h) \bigcap B(Q_i(x_{\ell}), 2 \theta) )
< C'_4 (n, k_2) e^{-a x_l},
\end{eqnarray}
then, we have
\begin{eqnarray}\label{eq:29-1}
\mbox{Vol}(G_3) < m^\rho \sum\limits_{\ell}  C'_4(n, k_2) e^{-a x_l} < 
C_4 (n, k_2) m^\rho e^{-a \eta}, 
\end{eqnarray}
since all these $P_i$ are sampled from the space $T \define \bigcup_{i=1}^{m^\rho} V_i$. Therefore, we can estimate the volume $G_3$ as:
\begin{eqnarray}\label{eq:29}
\mbox{Vol}(G_3) < C_4 (n, k_2) m^\rho e^{-a \eta} + \mbox{Vol}(\tilde{T}_h)
\end{eqnarray}

The total volume $T_h$ can be estiamted from $G$ via Eqs.~\eqref{eq:27},~\eqref{eq:210} and~\eqref{eq:29} as
\begin{eqnarray}\label{eq:211}
\mbox{Vol}(T_h) ~~<~~ \mbox{Vol}(G) &<& \mbox{Vol}( \tilde{T}_h ) + C_2(a, n, k_2) e^{-a(n-1)\eta} + C_3(n, k_2)e^{k_2 (n-1) \eta}  \nonumber \\
&  & + C_4 (n, k_2) m^\rho e^{-a \eta}.
\end{eqnarray}

By choosing $\eta$ as 
\begin{eqnarray}\label{eq:eta minimizer}
\eta = \frac{ \ln m^{\rho}}{ (n-1) (k_2 + a)},
\end{eqnarray}
we have the following bound between $\mbox{Vol}(T_h)$ and $\mbox{Vol}(T)$ by some computation
\begin{eqnarray}\label{eq:induction}
\mbox{Vol}(T_h) &< & \mbox{Vol}(T) + C_5(n, k_2)  \left( m^{\rho} \right)^{1 - \frac{a}{(n-1) (k_2 + a)}}. 
\end{eqnarray}

By finding proper $a$ according to Eqs.~\eqref{eq:25a} and~\eqref{eq:25b} (based on the value of $k_1$), from the volume induction relation provided by Eq.~\eqref{eq:induction}, this theorem is proved by setting 
\begin{eqnarray}
 \varpi(\rho, k_1,k_2,n)  = \frac{\rho a }{(n-1) (k_2 + a)}.
\end{eqnarray}
$\hfill \Box$

\subsection{Volume Estimation of $T$}\label{sec:Volume of T}

Because we have $T = \bigcup_{i=1}^{m^\rho} V_i$, where $0 \leq \rho \leq 1$, without loss of generality, we 
can assume the followin order of objects volume: $\mbox{Vol}(V_1) \geq \mbox{Vol}(V_2) \geq \cdots \geq  \mbox{Vol}(V_{m^\rho})$. We also assume that $\min\limits_{i \in \{1,2,\cdots,m\}} \mbox{Vol} (\partial V_i) \geq \beta$. Then, we have the following Lemma about the lower bound for the volume $T$. 

\begin{lemma}\label{lma:lower bound for vol T}
Let $M$ be a Hadamard manifold with dimension $n$ and sectional curvatures $K$ within in the range
$-k_2^2 \leq K \leq -k_1^2$, where $k_2 >  k_1 > 1$. We also have $V_1,\cdots,V_m$ are $m$ $\lambda$-convex 
sets in $M$ with $\lambda \leq k_2$ and assume that $\mbox{Vol}(V_1) \geq \mbox{Vol}(V_2) \geq \cdots \geq  \mbox{Vol}(V_{m^\rho})$. If $T = \bigcup_{i=1}^{m^\rho} V_i$, where $0 \leq \rho \leq 1$, the lower bound for the volume of $T$ can be derived as:
\begin{eqnarray}
\mbox{Vol}(T) \geq C_{lb}\frac{\lambda \beta }{k_2},
\end{eqnarray}
where $C_{lb}$ is a contant depending on the inradius of $V_1$ and $\beta$ is the lower bound for the surface area of $V_i$, where $i=1,2,\cdots,m$.
\end{lemma}
\textbf{Proof:}

Because $T=\bigcup_{i=1}^{m^\rho} V_i$, we have
\begin{eqnarray}
\mbox{Vol}(V_1) \leq \mbox{Vol}(T).
\end{eqnarray}
This Lemma is proved from Theorem 3 in~\cite{borisenko2002convex}.
$\hfill \Box$

\section{Covering Number Bounds by Volume}\label{sec:Covering Number Bounds by Volume}

In this section, we will try to upper bound the covering number of $T_h$ by the covering number of $T$ via covering number bounds by volume. We begin with the covering number definition. 
\begin{definition}
We use $N(A, d(\cdot, \cdot), \epsilon)$ to represent the covering number of the space $A$ with $\epsilon$-ball with distance metric function $d(\cdot, \cdot)$. Then, $N(A, d(\cdot, \cdot), \epsilon)$ can be defined as
\begin{eqnarray}
N(T, d(\cdot, \cdot), \epsilon) \define \min\{n: \mbox{the number of $\epsilon$-ball to cover $T$} \}.
\end{eqnarray}
\end{definition}

From Theorem 14.2 in~\cite{CoveringAndPacking2016}, we have
\begin{eqnarray}\label{eq:CoveringAndPacking2016}
\left(\frac{1}{\epsilon}\right)^n \frac{\mbox{Vol}(A) }{ \mbox{Vol}(\mathbf{B}) } \leq N(A, d(\cdot, \cdot), \epsilon) \leq \frac{\mbox{Vol}(A \oplus \frac{\epsilon}{2} \mathbf{B} ) }{ \mbox{Vol}( \frac{\epsilon}{2} \mathbf{B}) }
\leq_1 \left(\frac{3}{\epsilon}\right)^n \frac{\mbox{Vol}(A) }{ \mbox{Vol}(\mathbf{B}) }
\end{eqnarray}
where $\mathbf{B}$ is the unit norm ball, $\oplus$ is the Minkowski sum opertor, and the inequality $\leq_1$ is valid if $A$ is a convex set and $\epsilon \mathbf{B} \subset A$. 

We will have the following Lemma to upper bound the covering number of $T_h$ by the covering number of $T$ with respect to different types of $T$. We first define the following ratio between the volume upper bound for $T_h$ and the volume lower bound for $T$ as    
\begin{eqnarray}\label{eq:Ratio def}
R_{\mbox{Hada},n} \define \frac{  C_{ub} m^{1+\rho- \varpi(\rho, k_1,k_2,n)}    }{   C_{lb}\frac{\lambda \beta }{k_2}    },  
\end{eqnarray} 
from Thereom~\ref{thm:upper bound of vol Th} and Lemma~\ref{lma:lower bound for vol T}. 

\begin{lemma}\label{lma:Th convering num upper bound}
Let $M$ be a Hadamard manifold with dimension $n$ and sectional curvatures $K$ within in the range
$-k_2^2 \leq K \leq -k_1^2$, where $k_2 >  k_1 > 1$. We also have $V_1,\cdots,V_m$ are $m$ $\lambda$-convex 
sets in $M$ with $\lambda \leq k_2$ and $\bigcap_{i=1}^m V_i \neq \emptyset$. $\beta$ is the lower bound for the surface area of $V_i$, where $i=1,2,\cdots,m$. Suppose we select $m$ points, $P_1, P_2, \cdots, P_m$, such that all these $P_i$ for $1 \leq i \leq m$ are sampled from the space $T \define \bigcup_{i=1}^{m^\rho} V_i$, where $0 \leq \rho \leq 1$. Let us define $T_h$ as the convex hull of space $T$, then we have  
\begin{eqnarray}\label{eq1:lma:Th convering num upper bound}
N(T_h, d(\cdot, \cdot), \epsilon) &\leq&  R_{\mbox{Hada},n} 3^n  N(T, d(\cdot, \cdot), \epsilon),
\end{eqnarray}
where $R_{\mbox{Hada},n}$ is defined by Eq.~\eqref{eq:Ratio def}.
\end{lemma}
\textbf{Proof:}

Because $T_h$ is a convex set, we have
\begin{eqnarray}
N(T_h, d(\cdot, \cdot), \epsilon) &\leq&   \left(\frac{3}{\epsilon}\right)^n \frac{\mbox{Vol}(T_h) }{ \mbox{Vol}(\mathbf{B}) } \nonumber \\
&\leq_1&   \left(\frac{3}{\epsilon}\right)^n \frac{ R_{\mbox{Hada},n} \mbox{Vol}(T) }{ \mbox{Vol}(\mathbf{B}) } \nonumber \\
&=&  R_{\mbox{Hada},n} 3^n  \left(\frac{1}{\epsilon}\right)^n \frac{ \mbox{Vol}(T) }{ \mbox{Vol}(\mathbf{B}) } \nonumber \\
&\leq & R_{\mbox{Hada},n} 3^n   N(T, d(\cdot, \cdot), \epsilon),
\end{eqnarray}
where the equality $\leq_1$ comes from Thereom~\ref{thm:upper bound of vol Th} and Lemma~\ref{lma:lower bound for vol T}, and the first and the last 
inequalities are obtained from Eq.~\eqref{eq:CoveringAndPacking2016}. 
$\hfill \Box$

\section{Geometric Proof}\label{sec:Geometric Proof}

From Lemma~\ref{lma:Th convering num upper bound}, we are ready to prove Theorem 2.11.1 in~\cite{talagrand2021upper} geometrically for Hadamard manifold with dimension $n$. 

\begin{theorem}\label{thm:eq 2.130 in talagrand2021upper Book Poly}
Let $M$ be a Hadamard manifold with dimension $n$ and sectional curvatures $K$ within in the range
$-k_2^2 \leq K \leq -k_1^2$, where $k_2 >  k_1 > 1$. We also have $V_1,\cdots,V_m$ are $m$ $\lambda$-convex 
sets in $M$ with $\lambda \leq k_2$ and $\bigcap_{i=1}^m V_i \neq \emptyset$. $\beta$ is the lower bound for the surface area of $V_i$, where $i=1,2,\cdots,m$. Suppose we select $m$ points, $P_1, P_2, \cdots, P_m$, such that all these $P_i$ for $1 \leq i \leq m$ are sampled from the space $T \define \bigcup_{i=1}^{m^\rho} V_i$, where $0 \leq \rho \leq 1$. Let us define $T_h$ as the convex hull of space $T$, then we have  
\begin{eqnarray}
\gamma_\alpha(T_h, d(\cdot, \cdot)) \leq L_{\mbox{Hada}} \gamma_\alpha(T, d(\cdot, \cdot)),
\end{eqnarray}
where $L_{\mbox{Hada}}$ is the constant depending on the underlying geometry of the space $T$ and $\alpha$. 
\end{theorem}
\textbf{Proof:}

Because we have 
\begin{eqnarray}
\gamma_{\alpha}(T_h, d(\cdot, \cdot)) & \leqslant_{\alpha} & \int_0^{\infty}   \left(\log   N(T_h, d(\cdot, \cdot), \epsilon) \right)^{1/\alpha}   d \epsilon \nonumber \\
&\leq_1 &  \int_0^{\infty}   \left(\log  R_{\mbox{Hada},n} 3^n N(T, d(\cdot, \cdot), \epsilon) \right)^{1/\alpha}   d \epsilon \nonumber \\
&\leq &  \left(\frac{\log (R_{\mbox{Hada},n} 3^n)}{\log 2 }  + 1 \right)^{1/\alpha}\int_0^{\infty}   \left(\log  N(T, d(\cdot, \cdot), \epsilon) \right)^{1/\alpha}   d \epsilon \nonumber \\
& \leqslant_{\alpha}& \left(\frac{\log (R_{\mbox{Hada},n} 3^n)}{\log 2 }  + 1 \right)^{1/\alpha}\gamma_{\alpha}(T, d(\cdot, \cdot)), 
\end{eqnarray}
where the inequality $\leq_1$ comes from Lemma~\ref{lma:Th convering num upper bound}, and the first and last inequalities come from Theorem 1.2 in~\cite{talagrand2001majorizing}. This Lemma is proved. 
$\hfill \Box$

\bibliographystyle{IEEETran}
\bibliography{GeometricApproach_GC_Bib}

\end{document}